\documentclass{amsart}
\usepackage{graphicx}

\newtheorem{theorem}{Theorem}

\newtheorem{lemma}[theorem]{Lemma}

\newtheorem{prob}[theorem]{Problem}

\theoremstyle{remark}

\newtheorem{remark}[theorem]{Remark}

\newtheorem*{acknowledgement}{Acknowledgement}

\numberwithin{equation}{section}

\begin{document}

\title{Simplicial 2-spheres obtained from non-singular complete fans}

\author{Yusuke Suyama}
\address{Department of Mathematics, Graduate School of Science, Osaka City University,
3-3-138 Sugimoto, Sumiyoshi-ku, Osaka 558-8585 JAPAN}
\email{m13saU0r13@st.osaka-cu.ac.jp}



\date{\today}

\dedicatory{}

\begin{abstract}
We prove that a simplicial 2-sphere satisfying a certain condition
is the underlying simplicial complex of a 3-dimensional non-singular complete fan.
In particular, this implies that any simplicial 2-sphere with $\leq 18$ vertices
is the underlying simplicial complex of such a fan.
\end{abstract}

\maketitle

\section{Introduction}

A {\it rational strongly convex polyhedral cone} in $\mathbb{R}^n$
is a cone $\sigma$ spanned by finitely many vectors in $\mathbb{Z}^n$
which does not contain any non-zero linear subspace of $\mathbb{R}^n$.
A {\it fan} in $\mathbb{R}^n$ is a non-empty collection $\Delta$ of such cones
satisfying the following conditions:
\begin{enumerate}
\item If $\sigma \in \Delta$, then each face of $\sigma$ is in $\Delta$;
\item if $\sigma, \tau \in \Delta$, then $\sigma \cap \tau$ is a face of each.
\end{enumerate}
A fan $\Delta$ is {\it non-singular} if any cone in $\Delta$ is spanned by a part of a basis of $\mathbb{Z}^n$,
and {\it complete} if $\bigcup_{\sigma \in \Delta}\sigma=\mathbb{R}^n$.

A {\it toric variety} of complex dimension $n$ is a normal algebraic variety $X$ over $\mathbb{C}$
containing $(\mathbb{C}^*)^n$ as an open dense subset,
such that the natural action of $(\mathbb{C}^*)^n$ on itself extends to an action on $X$.
The category of toric varieties is equivalent to the category of fans (see \cite{O}).
A toric variety is smooth if and only if the corresponding fan is non-singular,
and compact if and only if the fan is complete.

Given a non-singular fan $\Delta$ with $m$ edges spanned by $v_1, \ldots, v_m \in \mathbb{Z}^n$,
we define its {\it underlying simplicial complex} as
\begin{equation*}
\{I \subset \{1, \ldots, m\} \mid \{v_i \mid i \in I\}\ {\rm spans\ a\ cone\ in}\ \Delta\}.
\end{equation*}
The underlying simplicial complex of an $n$-dimensional complete fan is a {\it simplicial $(n-1)$-sphere},
that is, a triangulation of the $(n-1)$-sphere.

For $n \geq 4$, a simplicial $(n-1)$-sphere is not always the underlying simplicial complex
of an $n$-dimensional non-singular complete fan (see \cite[Corollary 1.23]{DJ}).
On the other hand, successive equivariant blow-ups of $\mathbb{C}P^2$ produce
non-singular complete fans whose underlying simplicial complexes are all simplicial 1-spheres.
We consider the following problem:

\begin{prob}\label{problem}
Is any simplicial 2-sphere the underlying simplicial complex
of a 3-dimensional non-singular complete fan?
\end{prob}

No counterexamples to Problem \ref{problem} are currently known.
In this paper we give a partial affirmative answer to Problem \ref{problem}.
The {\it degree} of a vertex of a simplicial 2-sphere is the number of incident edges.

\begin{theorem}\label{theorem}
Let $K$ be a simplicial 2-sphere with $m_K$ vertices.
We denote the number of vertices of $K$ with degree $k$ by $p_K(k)$.
If $p_K(3)+p_K(4)+18 \geq m_K$,
then $K$ is the underlying simplicial complex of a 3-dimensional non-singular complete fan.
In particular, if $m_K \leq 18$, then $K$ is the underlying simplicial complex of such a fan.
\end{theorem}

The proof is done by reducing a given simplicial $2$-sphere
to another one in a collection of certain simplicial $2$-spheres with minimum degree $5$.
For each such simplicial 2-sphere,
we use a computer to find
a non-singular complete fan whose underlying simplicial complex is the simplicial 2-sphere.

The structure of the paper is as follows:
In Section 2, we give a complete list of the simplicial 2-spheres
with minimum degree 5 up to 18 vertices.
In Section 3, we prove Theorem \ref{theorem}.

\section{The simplicial 2-spheres with minimum degree 5 up to 18 vertices}

G. Brinkmann and B. D. McKay calculated the number of combinatorially different simplicial 2-spheres
with minimum degree 5 \cite{BM}:

\begin{table}[htbp]
\begin{center}
\begin{tabular}{|c||c|c|}
\hline
vertices & simplicial 2-spheres & simplicial 2-spheres with min. deg. 5 \\
\hline
 4 &           1 &  0 \\
\hline
 5 &           1 &  0 \\
\hline
 6 &           2 &  0 \\
\hline
 7 &           5 &  0 \\
\hline
 8 &          14 &  0 \\
\hline
 9 &          50 &  0 \\
\hline
10 &         233 &  0 \\
\hline
11 &       1,249 &  0 \\
\hline
12 &       7,595 &  1 \\
\hline
13 &      49,566 &  0 \\
\hline
14 &     339,722 &  1 \\
\hline
15 &   2,406,841 &  1 \\
\hline
16 &  17,490,241 &  3 \\
\hline
17 & 129,664,753 &  4 \\
\hline
18 & 977,526,957 & 12 \\
\hline
\end{tabular}
\caption{The number of simplicial 2-spheres.}
\label{number}
\end{center}
\end{table}

\begin{remark}
An $n$-dimensional {\it small cover} of a simple $n$-polytope
is a closed $n$-manifold $M$ with a locally standard $(\mathbb{Z}_2)^n$-action
such that the orbit space $M/(\mathbb{Z}_2)^n$ is the simple polytope.
It follows from Steinitz's theorem
that any simplicial 2-sphere is the boundary of a simplicial 3-polytope.
The dual of the simplicial 3-polytope is a simple 3-polytope $P$.
It follows from the four color theorem
that $P$ is the orbit space of a 3-dimensional small cover.
A 3-dimensional small cover of $P$ admits a hyperbolic structure
if and only if $P$ has no triangles or squares as facets,
that is, the original simplicial 2-sphere has no vertices with degree 3 or 4 \cite{DJ}.
Table \ref{number} shows that ``most" 3-dimensional small covers
do not admit any hyperbolic structure.
\end{remark}

We give a complete list of such simplicial 2-spheres
up to 18 vertices (see Tables \ref{table 1} and \ref{table 2}).
They are labeled as $\prod_{k \geq 5}k^{p(k)}$.
If there are more than one simplicial 2-spheres with the same label,
then we add (i), (ii), ... to the label.
Letters and $\star$ on vertices in Tables \ref{table 1} and \ref{table 2} are used in Section 3.

\begin{table}[htbp]
\begin{center}
\begin{tabular}{|c|c|c|}
\hline
& & \\
\includegraphics[width=3.2cm]{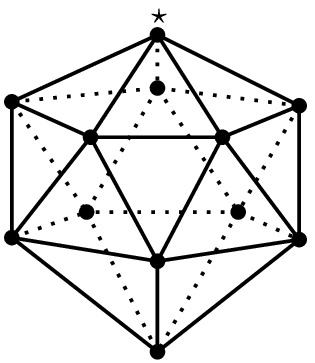} &
\includegraphics[width=3.2cm]{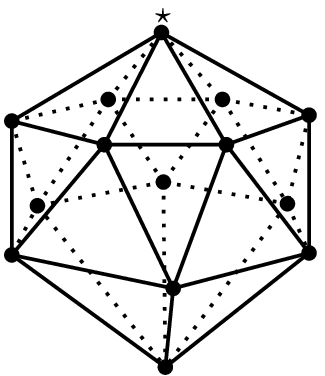} &
\includegraphics[width=3.2cm]{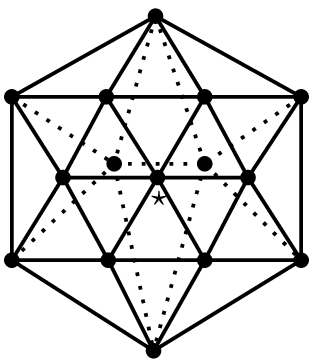} \\
$5^{12}$ & $5^{12}6^2$ & $5^{12}6^3$ \\
\hline
& & \\
\includegraphics[width=3.2cm]{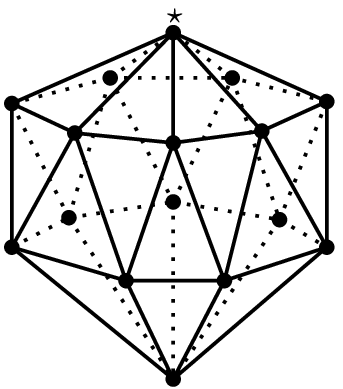} &
\includegraphics[width=3.2cm]{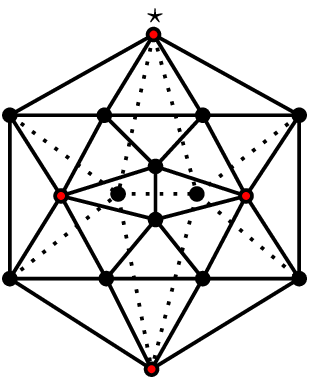} &
\includegraphics[width=3.2cm]{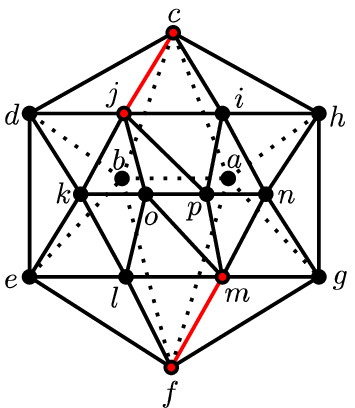} \\
$5^{14}7^2$ & $5^{12}6^4$ (i) & $5^{12}6^4$ (ii) \\
\hline
& & \\
\includegraphics[width=3.2cm]{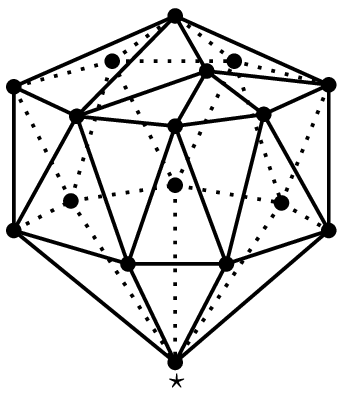} &
\includegraphics[width=3.2cm]{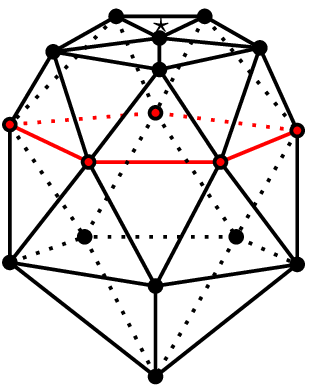} &
\includegraphics[width=3.2cm]{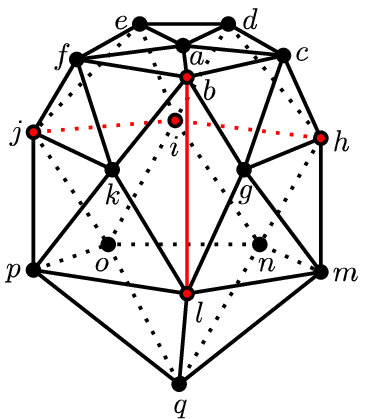} \\
$5^{13}6^37^1$ & $5^{12}6^5$ (i) & $5^{12}6^5$ (ii) \\
\hline
& & \\
\includegraphics[width=3.2cm]{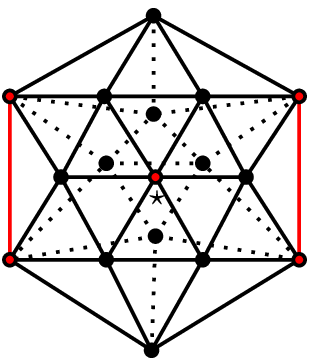} & & \\
$5^{12}6^5$ (iii) & & \\
\hline
\end{tabular}
\caption{The simplicial 2-spheres with minimum degree 5 up to 17 vertices.}
\label{table 1}
\end{center}
\end{table}

\begin{table}[htbp]
\begin{center}
\begin{tabular}{|c|c|c|}
\hline
& & \\
\includegraphics[width=3.2cm]{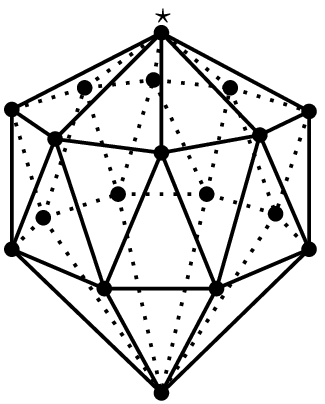} &
\includegraphics[width=3.2cm]{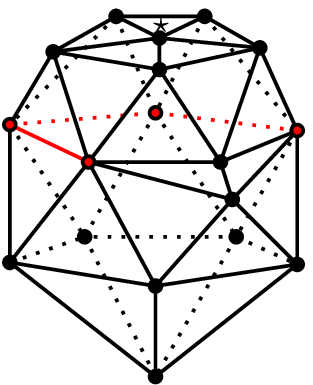} &
\includegraphics[width=3.2cm]{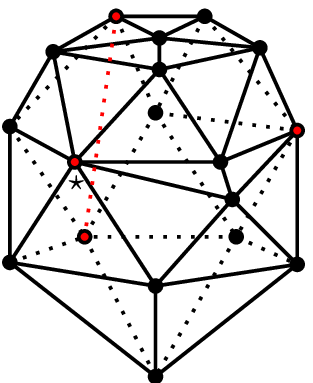} \\
$5^{16}8^2$ & $5^{14}6^27^2$ (i) & $5^{14}6^27^2$ (ii) \\
\hline
& & \\
\includegraphics[width=3.2cm]{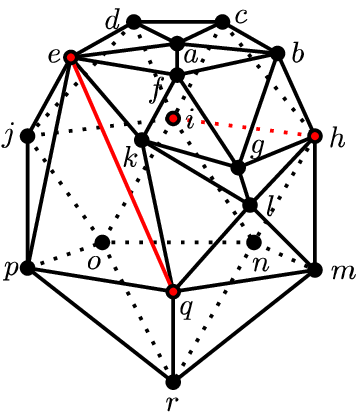} &
\includegraphics[width=3.2cm]{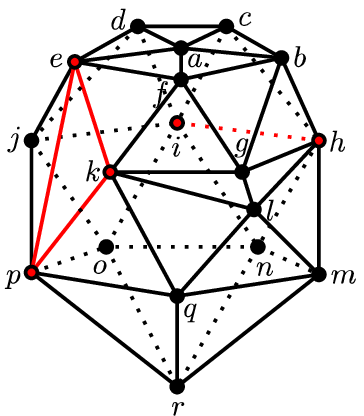} &
\includegraphics[width=3.2cm]{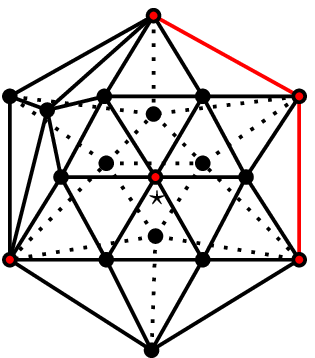} \\
$5^{14}6^27^2$ (iii) & $5^{13}6^47^1$ (i) & $5^{13}6^47^1$ (ii) \\
\hline
& & \\
\includegraphics[width=3.2cm]{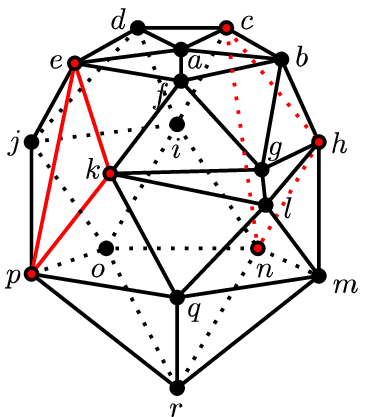} &
\includegraphics[width=3.2cm]{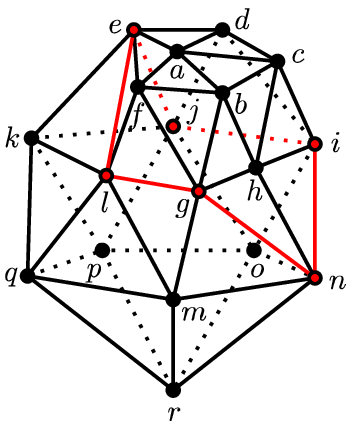} &
\includegraphics[width=3.2cm]{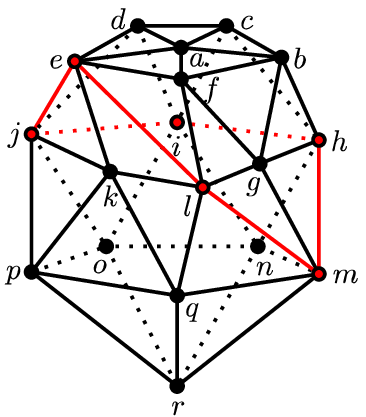} \\
$5^{12}6^6$ (i) & $5^{12}6^6$ (ii) & $5^{12}6^6$ (iii) \\
\hline
& & \\
\includegraphics[width=3.2cm]{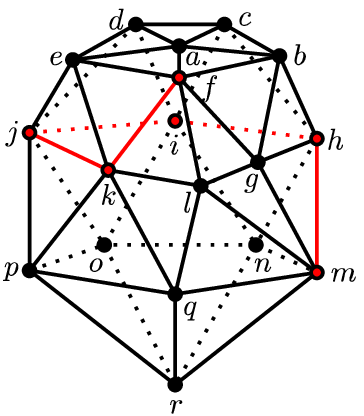} &
\includegraphics[width=3.2cm]{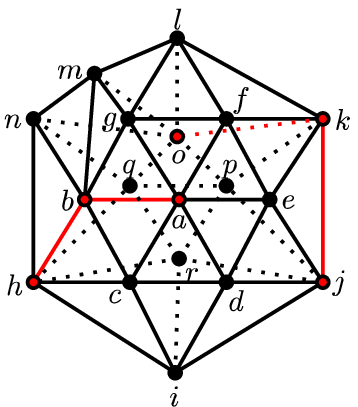} &
\includegraphics[width=3.2cm]{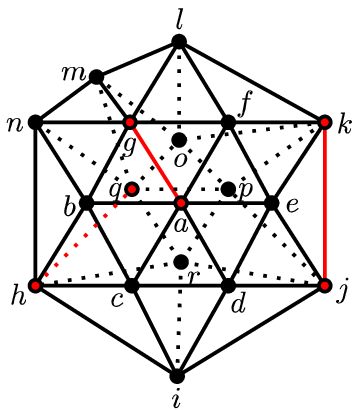} \\
$5^{12}6^6$ (iv) & $5^{12}6^6$ (v) & $5^{12}6^6$ (vi) \\
\hline
\end{tabular}
\caption{The simplicial 2-spheres with minimum degree 5 and 18 vertices.}
\label{table 2}
\end{center}
\end{table}

For each simplicial 2-sphere,
we consider the subcomplex consisting of the vertices with degree greater than or equal to 6
and the edges whose both endpoints have degree greater than or equal to 6
(red vertices and edges in Tables \ref{table 1} and \ref{table 2}).
These show that all simplicial 2-spheres in Tables \ref{table 1} and \ref{table 2}
are distinct except $5^{12}6^6$ (ii) and $5^{12}6^6$ (iii)
(they have the same subcomplex).

\begin{figure}[htbp]
\begin{center}
\includegraphics[width=3.5cm]{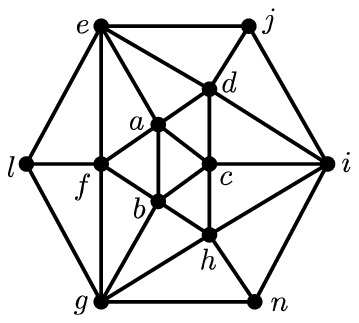}
\includegraphics[width=3.5cm]{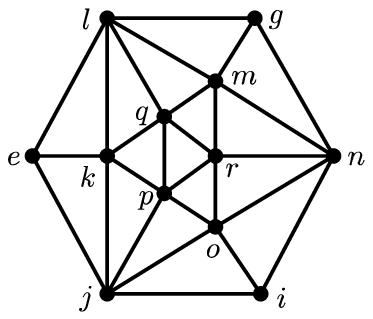}
\end{center}
\caption{Subcomplexes of $5^{12}6^6$ (ii).}
\label{sono12}
\end{figure}

\begin{figure}[htbp]
\begin{center}
\includegraphics[width=3.5cm]{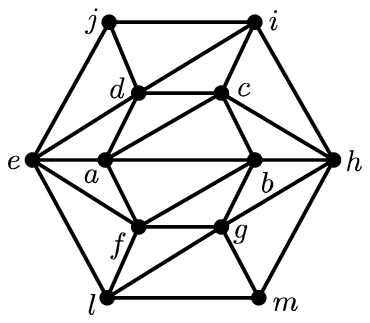}
\includegraphics[width=3.5cm]{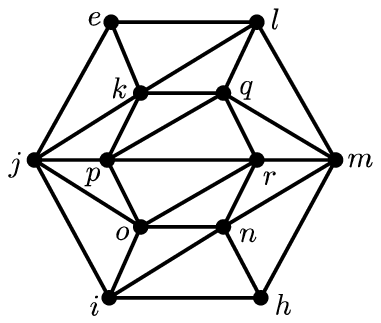}
\end{center}
\caption{Subcomplexes of $5^{12}6^6$ (iii).}
\label{sono34}
\end{figure}

Since the subcomplexes of $5^{12}6^6$ (ii) and $5^{12}6^6$ (iii) are cycles,
each cycle determines two subcomplexes surrounded by the cycle (see Figures \ref{sono12} and \ref{sono34}).
These are clearly distinct.

So all simplicial 2-spheres in Tables \ref{table 1} and \ref{table 2} are distinct.

For $m \leq 18$, the number of the simplicial 2-spheres with $m$ vertices
in Tables \ref{table 1} and \ref{table 2} agrees with the number in Table \ref{number}.
So this is a complete list of the simplicial 2-spheres
with minimum degree 5 up to 18 vertices.

\section{Proof of the Theorem \ref{theorem}}

Let $K$ be a simplicial 2-sphere with $m_K$ vertices.

\begin{figure}[htbp]
\begin{center}
\includegraphics[width=1.4cm]{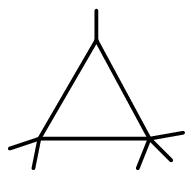}
$\stackrel{\rm (i)}{\longrightarrow}$
\includegraphics[width=1.4cm]{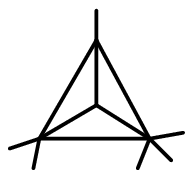}
\hspace{0.6cm}
\includegraphics[width=1.4cm]{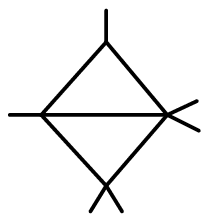}
$\stackrel{\rm (ii)}{\longrightarrow}$
\includegraphics[width=1.4cm]{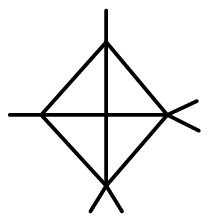}
\hspace{0.6cm}
\includegraphics[width=1.4cm]{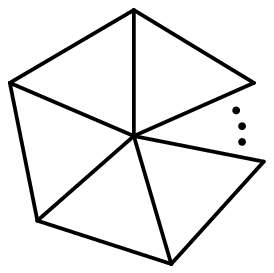}
$\stackrel{C_k}{\longrightarrow}$
\includegraphics[width=1.4cm]{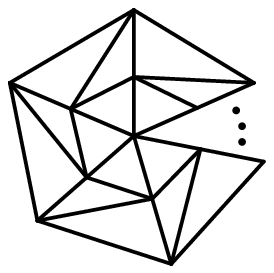} \\
For the operation $C_k$, the degree of the vertex in the center of the diagram is $k$.
\end{center}
\caption{Operations (i), (ii) and $C_k$.}
\label{operations}
\end{figure}

\begin{lemma}\label{lemma}
If $K$ is the underlying simplicial complex of a non-singular complete fan,
then a simplicial 2-sphere obtained from $K$ by an operation {\rm (i), (ii)} or $C_k \; (k \geq 5)$
is also the underlying simplicial complex of such a fan (see Figure \ref{operations}).
\end{lemma}

\begin{figure}[htbp]
\begin{center}
\includegraphics[height=4cm]{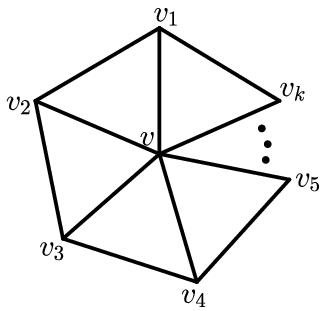}
$\stackrel{C_k}{\longrightarrow}$
\includegraphics[height=4cm]{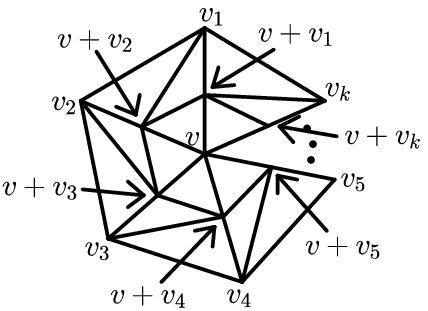}
\end{center}
\caption{An operation $C_k$.}
\label{C_k_3}
\end{figure}

\begin{proof}
Suppose that the three vertices of a 2-face of $K$ correspond to edge vectors $v_1, v_2, v_3 \in \mathbb{Z}^3$.
Then we have ${\rm det}(v_1, v_2, v_3)=1$.
We assign $v_1+v_2+v_3$ to the new vertex made by the operation (i).
The corresponding fan is non-singular and complete
since ${\rm det}(v_1, v_2, v_1+v_2+v_3)={\rm det}(v_2, v_3, v_1+v_2+v_3)={\rm det}(v_3, v_1, v_1+v_2+v_3)=1$.
Thus the lemma holds for an operation (i) (see Figure \ref{(i)_3}).

\begin{figure}[htbp]
\begin{center}
\includegraphics[height=2cm]{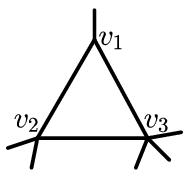}
$\stackrel{\rm (i)}{\longrightarrow}$
\includegraphics[height=2cm]{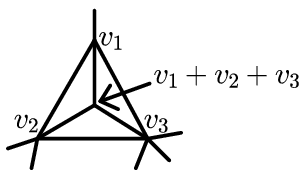}
\end{center}
\caption{An operation (i).}
\label{(i)_3}
\end{figure}

Suppose that $K$ contains a subcomplex in Figure \ref{(ii)_3}
and the vertices correspond to edge vectors $v_1, v_2, v_3, v_4 \in \mathbb{Z}^3$ as in Figure \ref{(ii)_3}.
Then we have ${\rm det}(v_1, v_2, v_3)={\rm det}(v_4, v_3, v_2)=1$.
We assign $v_2+v_3$ to the new vertex made by the operation (ii).
The corresponding fan is non-singular and complete
since ${\rm det}(v_1, v_2, v_2+v_3)={\rm det}(v_3, v_1, v_2+v_3)={\rm det}(v_2, v_4, v_2+v_3)
={\rm det}(v_4, v_3, v_2+v_3)=1$.
Thus the lemma holds for an operation (ii).

\begin{figure}[htbp]
\begin{center}
\includegraphics[height=2.5cm]{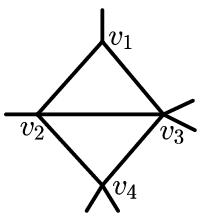}
$\stackrel{\rm (ii)}{\longrightarrow}$
\includegraphics[height=2.5cm]{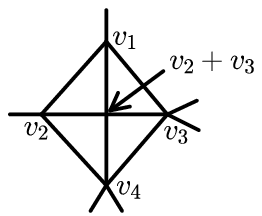}
\end{center}
\caption{An operation (ii).}
\label{(ii)_3}
\end{figure}

Suppose that $K$ contains a subcomplex in Figure \ref{C_k_3}
and the vertices correspond to edge vectors $v, v_1, \ldots, v_k \in \mathbb{Z}^3$ as in Figure \ref{C_k_3}.
Then we have ${\rm det}(v, v_i, v_{i+1})=1$ for any $i=1, \ldots, k$, where $v_{k+1}=v_1$.
For each $i=1, \ldots, k$, we assign $v+v_i$ to the new vertex between $v$ and $v_i$,
which is made by the operation $C_k$.
The corresponding fan is non-singular and complete
since ${\rm det}(v, v+v_i, v+v_{i+1})={\rm det}(v_i, v+v_{i+1}, v+v_i)={\rm det}(v_i, v_{i+1}, v+v_{i+1})=1$
for any $i=1, \ldots, k$.
Thus the lemma holds for an operation $C_k$.
This completes the proof.
\end{proof}

Now we prove Theorem \ref{theorem} by induction on $m_K$.
The tetrahedron is the only simplicial 2-sphere with 4 vertices,
which is the underlying simplicial complex of the fan of $\mathbb{C}P^3$.
Assume that $m_K \geq 5$.

(1) The case where there exists a vertex with degree 3.
All adjacent vertices have degree greater than or equal to 4,
since, if two vertices with degree 3 are adjacent, then $K$ must be the tetrahedron,
which contradicts $m_K \geq 5$.
Thus we can perform an inverse operation of (i) and we get a simplicial 2-sphere $K'$.
We see that $p_{K'}(3)+p_{K'}(4) \geq p_K(3)+p_K(4)-1$.
So we have $p_{K'}(3)+p_{K'}(4)+18 \geq p_K(3)+p_K(4)+18-1 \geq m_K-1=m_{K'}$.
$K'$ is the underlying simplicial complex of a non-singular complete fan by the induction hypothesis.
Hence $K$ is also the underlying simplicial complex of such a fan by Lemma \ref{lemma}.

(2) The case where there does not exist a vertex with degree 3
and there exists a vertex with degree 4.
Since all adjacent vertices have degree greater than or equal to 4,
we can perform an inverse operation of (ii) and we get a simplicial 2-sphere $K'$.
We see that $p_{K'}(3)+p_{K'}(4) \geq p_K(3)+p_K(4)-1$.
The same argument as (1) implies
that $K$ is the underlying simplicial complex of a non-singular complete fan.

(3) The case where there does not exist a vertex with degree 3 or 4.
The Euler relation implies that $\sum_{k \geq 3}(6-k)p_K(k)=12$ (see \cite[p.190]{O}).
This shows that $K$ must have a vertex with degree 5.
Since $m_K \leq p_K(3)+p_K(4)+18=18$ by assumption,
$K$ falls into 22 types in Tables \ref{table 1} and \ref{table 2}.

Suppose that $K$ has a vertex $v$ with degree $k \geq 5$ such that
any vertex adjacent to $v$ has degree 5, and
any vertex adjacent to a vertex adjacent to $v$ has degree greater than or equal to 5.
Then we can perform an inverse operation of $C_k$ and we get a simplicial 2-sphere $K'$.
Since $m_{K'}=m_K-k<18 \leq p_{K'}(3)+p_{K'}(4)+18$,
$K'$ is the underlying simplicial complex of a non-singular complete fan by the induction hypothesis.
Hence $K$ is also the underlying simplicial complex of such a fan by Lemma \ref{lemma}.

Each of $5^{12}$, $5^{12}6^5$ (i) and $5^{14}6^27^2$ (i) has such a vertex for $k=5$;
each of $5^{12}6^2$, $5^{12}6^3$, $5^{12}6^4$ (i), $5^{12}6^5$ (iii) and $5^{13}6^47^1$ (ii)
has such a vertex for $k=6$;
each of $5^{14}7^2$, $5^{13}6^37^1$ and $5^{14}6^27^2$ (ii) has such a vertex for $k=7$;
$5^{16}8^2$ has such a vertex for $k=8$
(these vertices are indicated by $\star$ in Tables \ref{table 1} and \ref{table 2}).
So they are the underlying simplicial complexes of non-singular complete fans.

We show that the rest of simplicial 2-spheres
$5^{12}6^4$ (ii), $5^{12}6^5$ (ii), $5^{14}6^27^2$ (iii), $5^{13}6^47^1$ (i) and $5^{12}6^6$ (i)--(vi)
are the underlying simplicial complexes of non-singular complete fans with a computer aid.
We assign vectors to the vertices as in Table \ref{vectors}.
They determine complete fans and it can be checked that all fans are non-singular by calculation.

\begin{table}[htbp]
\begin{center}
\begin{tabular}{|c||c|c|c|}
\hline
vertex & $5^{12}6^4$ (ii) & $5^{12}6^5$ (ii) & $5^{14}6^27^2$ (iii), $5^{13}6^47^1$ (i), $5^{12}6^6$ (i) \\
\hline
$a$ & $( 1, 0, 0)$ & $( 1, 0, 0)$ & $( 0,-1, 0)$ \\
\hline
$b$ & $( 0, 1, 0)$ & $( 1, 0, 1)$ & $( 1,-1, 0)$ \\
\hline
$c$ & $( 0, 0, 1)$ & $( 2,-1, 1)$ & $( 0,-1, 1)$ \\
\hline
$d$ & $(-1, 2,-1)$ & $( 3, 0,-1)$ & $(-1,-1, 1)$ \\
\hline
$e$ & $( 0,-1,-1)$ & $( 2, 1,-1)$ & $(-1,-1, 0)$ \\
\hline
$f$ & $( 1, 0,-1)$ & $( 1, 1, 0)$ & $(-1,-1,-1)$ \\
\hline
$g$ & $( 1,-1, 0)$ & $( 1,-1, 1)$ & $( 0,-1,-1)$ \\
\hline
$h$ & $( 1,-1, 1)$ & $( 2, 0,-1)$ & $( 1, 0, 0)$ \\
\hline
$i$ & $(-1, 0, 1)$ & $( 1, 1,-1)$ & $( 0, 0, 1)$ \\
\hline
$j$ & $(-1, 1, 0)$ & $( 0, 1, 0)$ & $(-1, 0, 1)$ \\
\hline
$k$ & $(-1, 1,-1)$ & $( 0, 0, 1)$ & $(-1, 0,-1)$ \\
\hline
$l$ & $( 0,-2,-1)$ & $( 0,-1, 1)$ & $( 0, 0,-1)$ \\
\hline
$m$ & $( 1,-1,-1)$ & $( 2,-1, 0)$ & $( 0, 1,-1)$ \\
\hline
$n$ & $( 0,-1, 1)$ & $( 1, 0,-1)$ & $( 1, 1, 0)$ \\
\hline
$o$ & $( 0,-1, 0)$ & $( 0, 1,-1)$ & $( 0, 1, 1)$ \\
\hline
$p$ & $( 0,-2, 1)$ & $(-1, 1, 0)$ & $(-1, 0, 0)$ \\
\hline
$q$ &              & $(-1, 0, 0)$ & $(-1, 1,-1)$ \\
\hline
$r$ &              &              & $( 0, 1, 0)$ \\
\hline
\end{tabular}
\begin{tabular}{|c||c|c|c|c|c|}
\hline
vertex & $5^{12} 6^6$ (ii) & $5^{12} 6^6$ (iii) & $5^{12} 6^6$ (iv) & $5^{12} 6^6$ (v) & $5^{12} 6^6$ (vi) \\
\hline
$a$ & $( 1, 0, 0)$ & $( 1, 0, 0)$ & $( 1, 0, 0)$ & $( 0,-1, 0)$ & $( 0,-1, 0)$ \\
\hline
$b$ & $( 3, 0,-1)$ & $( 3, 0,-1)$ & $( 3, 0,-1)$ & $(-1, 1,-1)$ & $(-1, 0,-1)$ \\
\hline
$c$ & $( 2, 1,-1)$ & $( 2, 1,-1)$ & $( 2, 1,-1)$ & $( 0,-2,-1)$ & $( 0,-2,-1)$ \\
\hline
$d$ & $( 1, 1, 0)$ & $( 1, 1, 0)$ & $( 1, 1, 0)$ & $( 1,-1,-1)$ & $( 1,-1,-1)$ \\
\hline
$e$ & $( 3, 0, 1)$ & $( 1, 0, 1)$ & $( 1, 0, 1)$ & $( 0,-1, 1)$ & $( 0,-1, 1)$ \\
\hline
$f$ & $( 3,-1, 1)$ & $( 3,-1, 1)$ & $( 2,-1, 1)$ & $(-1, 0, 1)$ & $(-1, 0, 1)$ \\
\hline
$g$ & $( 2, 0,-1)$ & $( 2, 0,-1)$ & $( 2, 0,-1)$ & $(-1, 1, 0)$ & $(-1, 1, 0)$ \\
\hline
$h$ & $( 1, 1,-1)$ & $( 1, 1,-1)$ & $( 1, 1,-1)$ & $( 0,-1,-1)$ & $( 0,-1,-1)$ \\
\hline
$i$ & $( 0, 1, 0)$ & $( 0, 1, 0)$ & $( 0, 1, 0)$ & $( 1, 0,-1)$ & $( 1, 0,-1)$ \\
\hline
$j$ & $( 1, 0, 1)$ & $( 0, 0, 1)$ & $( 0, 0, 1)$ & $( 1,-1, 0)$ & $( 1,-1, 0)$ \\
\hline
$k$ & $( 1,-1, 1)$ & $( 1,-1, 1)$ & $( 1,-1, 1)$ & $( 1,-1, 1)$ & $( 1,-1, 1)$ \\
\hline
$l$ & $( 2,-1, 1)$ & $( 2,-1, 1)$ & $( 3,-1, 0)$ & $( 0, 0, 1)$ & $( 0, 0, 1)$ \\
\hline
$m$ & $( 1, 0,-1)$ & $( 1, 0,-1)$ & $( 1, 0,-1)$ & $(-1, 2, 0)$ & $(-1, 2, 2)$ \\
\hline
$n$ & $(-1, 1, 0)$ & $( 0, 1,-1)$ & $( 0, 1,-1)$ & $(-1, 2,-1)$ & $(-2, 2,-1)$ \\
\hline
$o$ & $( 0, 0, 1)$ & $(-1, 1, 0)$ & $(-1, 1, 0)$ & $( 0, 1, 2)$ & $( 0, 1, 2)$ \\
\hline
$p$ & $( 0,-1, 1)$ & $( 0,-1, 1)$ & $( 0,-1, 1)$ & $( 0, 1, 1)$ & $( 0, 1, 1)$ \\
\hline
$q$ & $( 2,-1, 0)$ & $( 2,-1, 0)$ & $( 2,-1, 0)$ & $(-1, 2,-2)$ & $(-1, 1,-1)$ \\
\hline
$r$ & $(-1, 0, 0)$ & $(-1, 0, 0)$ & $(-1, 0, 0)$ & $( 0, 1, 0)$ & $( 0, 1, 0)$ \\
\hline
\end{tabular}
\caption{Assigning vectors to the vertices.}
\label{vectors}
\end{center}
\end{table}

\begin{figure}[htbp]
\begin{center}
\includegraphics[width=8.5cm]{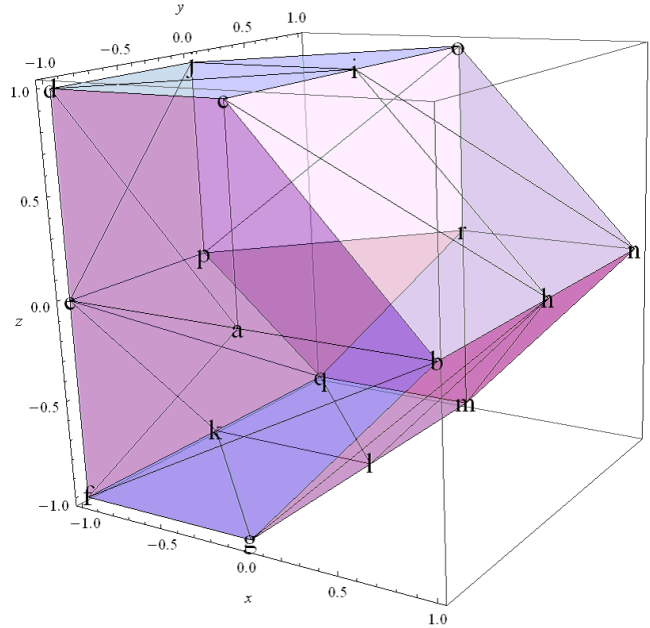}
\end{center}
\caption{$5^{14}6^27^2$ (iii).}
\label{mathematica}
\end{figure}

For example, we show that $5^{14}6^27^2$ (iii) is the underlying simplicial complex of a non-singular complete fan.
Vectors in Table \ref{vectors} determine a 3-dimensional complete fan.
Its underlying simplicial complex is illustrated in Figure \ref{mathematica},
which confirms that there are no overlaps among the 3-dimensional cones.
Calculating determinants, say $\det(a, b, c)=1$, we see that every cone is non-singular.

\begin{acknowledgement}
The author wishes to thank Professor Mikiya Masuda
for his valuable advice and continuing support,
and Professors Hiroshi Sato, Tadao Oda and Masanori Ishida
for useful comments on Problem \ref{problem} in the introduction.
\end{acknowledgement}

\end{document}